\documentclass[10pt]{amsart}

\usepackage{enumerate}
\usepackage{graphics}
\usepackage{graphicx}
\usepackage{amsmath}
\usepackage{amssymb}
\usepackage{amstext}
\usepackage{times}
\usepackage[latin1]{inputenc}
\usepackage{latexsym}
\usepackage{amsthm}
\usepackage{amscd}

\def\R{\mathbb{R}}
\def\L{\mathbb{L}}
\def\C{\mathbb{C}}
\def\N{\mathbb{N}}

\def\H{\mathbb{H}}
\def\B{{\tt B}}
\def\b{{\tt b}}
\def\esf{\mathbb{S}}

\newcommand{\Real}{\mbox{\rm Re}}

\newcommand{\dist}{\operatorname{dist}}
\newcommand{\length}{\operatorname{length}}
\newcommand{\diam}{\operatorname{diam}}

\newcommand{\intc}{\operatorname{Int}}
\newcommand{\escpro}[1]{\langle {#1} \rangle}
\newcommand{\norma}[1]{\| {#1} \|}
\newcommand{\escproR}[1]{\langle {#1} \rangle_{0}}
\newcommand{\normaR}[1]{\| {#1} \|_{0}}

\def\a{{\alpha}}
\def\om{{\omega}}

\def\g{{\gamma}}
\def\l{{\lambda}}
\def\de{{\delta}}

\def\Om{{\Omega}}

\def\ep{{\epsilon}}
\def\ve{{\varepsilon}}
\def\ri{{\rm i}}
\def\NN{{\mathcal{N}}}
\def\cte{{\tt const}}

\newtheorem{lemma}{Lemma}
\newtheorem{remark}{Remark}
\newtheorem{theorem}{Theorem}

\newtheorem{definition}{Definition}
\newtheorem{claim}{Claim}

\begin{document}

\title{On the Calabi-Yau problem for maximal surfaces in $\L^3$}
\author{Antonio Alarcón}

\thanks{This research is partially supported by MEC-FEDER Grant no. MTM2004 - 00160.}
\date{\today}

\address{Departamento de Geometría y Topología
\hfill\break\indent Universidad de Granada, \hfill\break\indent
18071, Granada \hfill\break\indent Spain}

\email{alarcon@ugr.es}

\begin{abstract}
In this paper we construct an example of a maximal surface in the
Lorentz-Minkowski space $\L^3,$ which is bounded by a hyperboloid
and weakly complete in the sense explained by Umehara and Yamada
\cite{Um-Ya}.
\\

\noindent {\em 2000 Mathematics Subject Classification:} Primary 53C50; Secondary 53C42, 53A10, 53B30.

\noindent {\em Keywords:} Complete maximal immersions, maximal
surfaces with singularities.
\end{abstract}

\maketitle



\section{Introduction}

A maximal hypersurface in a Lorentzian manifold is a spacelike
hypersurface with zero mean curvature. Besides of their
mathematical interest these hypersurfaces and more generally those
having constant mean curvature have a significant importance in
physics \cite{Kiehn,Kiehn2,MT}. When the ambient space is the
Minkowski space $\L^n,$ one of the most important results is the
proof of a Bernstein-type theorem for maximal hypersurfaces in
$\L^n$. Calabi \cite{Calabi} proved that the only complete
hypersurfaces with zero mean curvature in $\L^3$ (i.e. maximal
surfaces) and $\L^4$ are spacelike hyperplanes, solving the so
called Bernstein-type problem in dimensions 3 and 4. Cheng and Yau
\cite{Ch-Yau} extended this result to $\L^n,$ $n \geq 5.$ It is
therefore meaningless to consider global problems on maximal and
everywhere regular hypersurfaces in $\L^n.$ In contrast, there
exists a lot of results about existence of non-flat maximal
surfaces with singularities \cite{Isa-PL1,Isa-PL-RS,LLS}.
\\

It is well known the close relationship between maximal surfaces
in $\L^3$ and minimal surfaces in $\R^3$ (see Remark \ref{minmax}
in page \pageref{minmax}). This fact let us solve some problems on
maximal surfaces by solving the analogous ones for minimal
surfaces, and vice versa. This is not the case of the Calabi-Yau
problem. In 1965 Calabi asked whether or not it is possible for a
complete minimal surface in $\R^3$ to be bounded. Much work has
been done on it over the past four decades. The most important
result in this line was obtained by Nadirashvili \cite{nadi}, who
constructed a complete minimal surface in the unit ball of $\R^3.$
See \cite{salvador} for more information about this topic. From a
Nadirashvili's surface and using the relationship between maximal
and minimal surfaces, we can obtain as most the existence of a
weakly complete maximal surface contained in a cylinder of $\L^3.$
Here, we use the concept of weakly completeness (see definition
\ref{def: weakly} in page \pageref{def: weakly}) that was
introduced by Umehara and Yamada \cite{Um-Ya}.
\\

In this paper, we construct an example of a weakly complete
maximal surface in $\L^3$ with singularities, which is bounded by
a hyperboloid. We would like to point out that our example does
not have branch points, all the singularities are of lightlike
type (see definition \ref{def: light} in page \pageref{def:
light}).

More precisely, we prove the following existence theorem.

\begin{theorem}\label{teorema}
There exists a weakly complete conformal maximal immersion with
lightlike singularities of the unit disk into the set
$\{(x,y,z)\in\L^3\;|\;x^2+y^2-z^2<-1\}.$
\end{theorem}

For several reasons, lightlike singularities of maximal surfaces
in $\L^3$ are specially interesting. This kind of singularities
are more attractive than branch points, in the sense that they
have a physical interpretation \cite{Kiehn,Kiehn2}. At these
points, the limit tangent plane is lightlike, the curvature blows
up and the Gauss map has no well defined limit. However, as in the
case of minimal surfaces, if we allow branch points, then proving
the analogous result of Theorem \ref{teorema} has less technical
difficulties.

The fundamental tools used in the proof of this result (Runge's
theorem and the López-Ros transformation) are those that
Nadirashvili utilized to construct the first example of a complete
bounded minimal surface in $\R^3.$ Improvements of his technique
have generated a lot of literature on the Calabi-Yau problem for
minimal surfaces in $\R^3$ \cite{MM-Convex2, MN, AFM}.
\\

Similarly to the case of minimal surfaces, it would be stimulating
to look for an additional property for a weakly complete bounded
maximal surface: properness. In order to achieve it, the technique
showed in this paper could be combined with the reasonings used in
the construction \cite{dga} of a proper conformal maximal disk in
$\L^3,$ following the ideas of \cite{MM-Convex}. The main
objection of this argument is that the best result known about the
convex hull property for maximal surfaces \cite{leo-rosa} needs
the control of the image of the singularities of the surface. This
problem will be studied in \cite{yo-isa}.


\section{Background and notation}

\subsection{The Lorentz-Minkowski three space}

We denote by $\L^3$ the three dimensional Lorentz-Minkowski space
$(\R^3,\escpro{\cdot,\cdot}),$ where
$\escpro{\cdot,\cdot}=dx_1^2+dx_2^2-dx_3^2.$ The Lorentzian {\em
norm} is given by $\|(x_1,x_2,x_3)\|^2=x_1^2+x_2^2-x_3^2,$ and
$\|x\|=\rm{sign}(\|x\|^2)\sqrt{|\|x\|^2|}.$ We say that a vector
$v\in\R^3\setminus \{(0,0,0)\}$ is spacelike, timelike or
lightlike if $\|v\|^2$ is positive, negative or zero,
respectively. The vector $(0,0,0)$ is spacelike by definition. A
plane in $\L^3$ is spacelike, timelike or lightlike if the induced
metric is Riemannian, non degenerate and indefinite or degenerate,
respectively.

In order to differentiate between $\L^3$ and $\R^3$, we denote
$\R^3=(\R^3,\escproR{\cdot,\cdot}),$ where $\escproR{\cdot,\cdot}$
is the usual metric of $\R^3$, i.e., $\escproR{\cdot,\cdot}=dx_1^2+dx_2^2+dx_3^2.$
We also denote the Euclidean norm by $\normaR{\cdot}.$
\\

By an (ordered) $\L^3$-orthonormal basis we mean a basis of $\R^3,$ $\{u,v,w\},$ satisfying
\begin{itemize}
\item $\escpro{u,v}=\escpro{u,w}=\escpro{v,w}=0$;
\item $\norma{u}=\norma{v}=-\norma{w}=1$.
\end{itemize}
Notice that $u$ and $v$ are spacelike vectors whereas $w$ is
timelike.
\\

We call $\H^2:=\{(x_1,x_2,x_3)\in\R^3 \;|\;
x_1^2+x_2^2-x_3^2=-1\}$ the hyperbolic sphere in $\L^3$ of
constant intrinsic curvature $-1.$ Notice that $\H^2$ has two
connected components $\H^2_+:=\H^2\cap\{x_3\geq 1\}$ and
$\H^2_-:=\H^2\cap \{x_3\leq -1\}.$ The stereographic projection
$\eta$ for $\H^2$ from the point $(0,0,1)\in\H^2_+$ is the map
$\eta:\H^2\to \C\cup\{\infty\}\setminus \{|z|=1\}$ given by
\[
\eta (x_1,x_2,x_3)=\frac{x_1+\ri x_2}{1-x_3} \;,\quad
\eta(0,0,1)=\infty\;.
\]
Notice that $\eta(\H^2_+)=\{|z|>1\}$ and $\eta(\H^2_-)=\{|z|<1\}.$

Given $r\geq 0,$ we denote by $\B(r)$ as the lower convex domain
determined by the set $\{\|x\|=-r\},$ i.e.,
\[
\B(r)=\{(x_1,x_2,x_3)\in\R^3 \;|\; \norma{(x_1,x_2,x_3)}<-r\;,\;
x_3<-r \}\;.
\]
We also denote $\b(r)=\partial \B(r).$ Observe that
$\b(1)=\H^2_-.$ Moreover, if $r_1<r_2,$ then
$\overline{\B(r_2)}\subset \B(r_1)$ and $\b(r_1)\cap
\b(r_2)=\emptyset.$
\\

Finally, we define the maps $\NN:\B(0)\to\H_+^2$ and
$\NN_0:\B(0)\to\esf^2$ by the following way. Consider $p\in B(0)$
and label $r=\|p\|<0.$ Let $\NN^r:\b(r)\to\H^2_+$ and
$\NN_0^r:\b(r)\to \esf^2$ be the outward pointing $\L^3$-normal
Gauss map and the Euclidean outward pointing unit normal of
$b(r),$ respectively. Then, we define
\[
\NN(p)=\NN^r(p)\;,\quad \NN_0(p)=\NN_0^r(p) \;.
\]\label{pag: normales}
Equivalently, $\NN(p)=-p/\norma{p}$ and
$\NN_0(p)=\mathcal{J}(p)/\normaR{p},$ where
$\mathcal{J}(p_1,p_2,p_3)=(p_1,p_2,-p_3).$ Hence, both maps are
differentiable and
$\NN_0(p)=-\mathcal{J}(\NN(p))/\normaR{\NN(p)}.$

\subsection{Maximal surfaces}

Any conformal maximal immersion $X:M\to\L^3$ is given by a triple
$\Phi=(\Phi_1,\Phi_2,\Phi_3)$ of holomorphic 1-forms defined on
the Riemann surface $M,$ having no common zeros and satisfying
\begin{equation}\label{2}
|\Phi_1|^2+|\Phi_2|^2-|\Phi_3|^2\neq 0\;;
\end{equation}
\begin{equation}\label{conforme}
\Phi_1^2+\Phi_2^2-\Phi_3^2=0\;;
\end{equation}
and all periods of the $\Phi_j$ are purely imaginary. Here we
consider $\Phi_i$ to be a holomorphic function times $dz$ in a
local parameter $z.$ Then, the maximal immersion $X:M\to\L^3$ can
be parameterized by $z\mapsto \Real \int^z \Phi.$ The above triple
is called the Weierstrass representation of the maximal immersion
$X.$ Usually, the second requirement \eqref{conforme} is
guaranteed by the introduction of the formulas
\[
\Phi_1=\frac{\ri}{2}(1-g^2)\eta\;,\quad \Phi_2=-\frac12(1+g^2)\eta\;,\quad \Phi_3=g\eta
\]
for a meromorphic function $g$ with $|g(p)|\neq 1,$ $\forall p\in
M,$ (the stereographically projected Gauss map) and a holomorphic
1-form $\eta.$ We also call $(g,\eta)$ or $(g,\Phi_3)$ the
Weierstrass representation of $X.$

\begin{remark}\label{minmax}
If $(\Phi_1,\Phi_2,\Phi_3)$ is the Weierstrass representation of a
maximal surface, then $(\ri \Phi_1,\ri \Phi_2, \Phi_3)$ are the
Weierstrass data of a minimal surface in $\R^3$ \cite{osserman}.
Moreover, both surfaces have the same meromorphic Gauss map $g.$
\end{remark}

We are going to deal with maximal immersions with lightlike
singularities, according with the following definition.
\begin{definition}\label{def: light}
A point $p\in M$ is a lightlike singularity of the immersion $X$
if it is not a branch point and $|g(p)|=1.$
\end{definition}

In this article, all the maximal immersions are defined on simply
connected domains of $\C,$ thus the Weierstrass 1-forms have no
periods and so the only requirements are \eqref{2} at the points
that are not singularities, and \eqref{conforme}. In this case,
the differential $\eta$ can be written as $\eta =f(z)dz.$ The
metric of $X$ can be expressed as
\begin{equation}\label{metrica}
ds^2=\frac12 (|\Phi_1|^2+|\Phi_2|^2-|\Phi_3|^2)=\big( \frac12 (1-|g|^2)|f||dz| \big)^2\;.
\end{equation}
We use a concept of completeness that is less exigent than the
classical one. The following definition was given by Umehara and
Yamada \cite{Um-Ya}.

\begin{definition}\label{def: weakly}
A maximal immersion $X:M\to\L^3$ is weakly complete if the Riemann
surface $M$ is complete with the metric
\begin{equation}\label{metricaUY}
d\sigma^2=\frac12 (|\Phi_1|^2+|\Phi_2|^2+|\Phi_3|^2)=\big( \frac12
(1+|g|^2)|f||dz| \big)^2\;.
\end{equation}
The metric $d\sigma^2$ will be called the lift metric of $X.$
\end{definition}

The Euclidean metric on $\C$ is denoted as $\escpro{,}=|dz|^2.$
Note that $ds^2=(\l_X)^2\,|dz|^2$ and
$d\sigma^2=(\l^0_X)^2\,|dz|^2$ where the conformal coefficients
$\l_X$ and $\l^0_X$ are given by \eqref{metrica} and
\eqref{metricaUY}, respectively.

\begin{remark}\label{rem: branch}
Observe that if $X$ has a singularity of lightlike type in a point
$z\in M,$ then $\l_X(z)=0$ but $\l_X^0(z)\neq 0.$ On the other
hand, if $z$ is a branch point of $X,$ one has
$\l_X(z)=0=\l_X^0(z).$
\end{remark}

Along this paper, we use some $\L^3$-orthonormal bases. Given
$X:\Omega\to\L^3$ a maximal immersion and $S$ an
$\L^3$-orthonormal basis, we write the Weierstrass data of $X$ in
the basis $S$ as
\[
\Phi_{(X,S)}=(\Phi_{(1,S)},\Phi_{(2,S)},\Phi_{(3,S)})\;, \quad
f_{(X,S)}\;, \quad g_{(X,S)}\;, \quad \eta_{(X,S)}\;.
\]

In the same way, given $v\in\R^3,$ we denote by $v_{(k,S)}$ the
$k$th coordinate of $v$ in $S.$ We also represent by
$v_{(*,S)}=(v_{(1,S)},v_{(2,S)})$ the first two coordinates of $v$
in the basis $S.$
\\

Given a curve $\a$ in $\Om,$ by $\length (\a,ds)$ we mean the
length of $\a$ with respect to the metric $ds.$ Let $W\subset \Om$
be a subset, then we define
\begin{itemize}
\item $\dist_{(W,ds)}(p,q)=\inf \{\length (\a,ds) \;|\; \a:[0,1]\to W,$ $\a(0)=p,$ $\a(1)=q\},$ for any $p,q\in W.$

\item $\dist_{(W,ds)}(U,V)= \inf \{ \dist_{(W,ds)}(p,q) \;|\; p\in U,$ $q\in V\},$ for any $U,V\subset W.$
\end{itemize}

Given a domain $D\subset \C,$ we say that a function, or a 1-form,
is harmonic, holomorphic, meromophic, ... on $\overline{D}$, if it
is harmonic, holomorphic, meromorphic, ... on a domain containing
$\overline{D}.$

Let $P$ be a simple closed polygonal curve in $\C$. By $\intc P$
we mean the bounded connected component of $\C\setminus P.$ For a
small enough $\xi>0,$ we denote by $P^\xi$ as the parallel
polygonal curve in $\intc P,$ satisfying that the distance between
parallel sides is equal to $\xi.$ Whenever we write $P^\xi$ we are
assuming that $\xi$ is small enough to define the polygon
properly.

\subsection{The López-Ros transformation}

The proof of Lemma \ref{lema} exploits what has come to be call
the López-Ros transformation. If $(g,f)$ are the Weierstrass data
of a maximal immersion $X:\Omega\to\L^3$ (being $\Omega$ simply
connected), we define on $\Omega$ the data
\[
\widetilde{g}=\frac{g}{h}\;,\quad \widetilde{f}=f\,h\;,
\]
where $h:\Omega\to\C$ is a holomorphic function without zeros.
Observe that the new meromorphic data satisfy \eqref{2} at the
regular points, and \eqref{conforme}, so the new data define a
maximal immersion (possibly with different lightlike
singularities) $\widetilde{X}:\Omega\to\L^3.$ This method provides
us with a powerful and natural tool for deforming maximal
surfaces. One of the most interesting properties of the resulting
surface is that the third coordinate function is preserved.


\section{Proof of Theorem \ref{teorema}}

In order to prove Theorem \ref{teorema} we will apply the
following technical Lemma. It will be proved later in Section
\ref{sec: lema}.

\begin{lemma}\label{lema}
Consider $r>0,$ $P$ a polygon in $\C$ and $X:\overline{\intc P}\to
\L^3$ a conformal maximal immersion (possibly with lightlike
singularities) satisfying
\begin{equation}\label{ecu: X-cota}
X(\overline{\intc P})\subset \B(r)\;.
\end{equation}
Let $\ep$ and
$s$ be positive constants with $\sqrt{r^2-4s^2}-\ep>0.$ Then,
there exist a polygon $Q$ and a conformal maximal immersion
(possibly with lightlike singularities) $Y:\overline{\intc
Q}\to\L^3$ such that
\begin{enumerate}[\rm ({L}.1)]
\item $\overline{\intc P^\ep}\subset \intc Q\subset \overline{\intc Q}\subset \intc P.$

\item $s<\dist_{(\overline{\intc Q},d\sigma_Y^2)}(P^\ep,Q),$ where
$d\sigma_Y^2$ is the lift metric associated to the immersion $Y.$

\item $Y(\overline{\intc Q})\subset \B(R),$ where
$R=\sqrt{r^2-4s^2}-\ep.$

\item $\|Y-X\|_0<\ep$ in $\overline{\intc P^\ep}.$
\end{enumerate}
\end{lemma}

Using this Lemma, we construct a sequence of immersions
$\{\psi_n\}_{n\in\N}$ that converges to an immersion $\psi$ which
proves Theorem \ref{teorema}, up to a reparametrization of its
domain.

First of all, we consider a sequence of reals $\{\a_n\}_{n\in\N}$
satisfying
\[
\prod_{k=1}^\infty \a_k=\frac12\;,\quad 0<\a_k<1\;,\quad \forall
k\in\N\;.
\]
Moreover, we choose $r_1>1$ large enough so that the sequence
$\{r'_n\}_{n\in\N}$ given by
\[
r'_1=r_1\;,\quad r'_n=\sqrt{(r'_{n-1})^2-(2/n)^2}-\frac1{n^2}
\]
satisfies
\begin{equation}\label{ecu: r'}
r'_n>1\;,\quad \forall n\in\N\;.
\end{equation}

Now, we are going to construct a sequence
$\{\Upsilon_n\}_{n\in\N},$ where the element
\[
\Upsilon_n=\{P_n,\psi_n,\ep_n,\xi_n\}
\]
is composed of a polygon $P_n,$ a conformal maximal immersion
$\psi_n:\overline{\intc P_n}\to\L^3,$ and $\ep_n<\frac1{n^2},$ and
$\xi_n$ are positive real numbers. We will choose $\ep_n$ and
$\xi_n$ so that the sequences $\{\ep_n\}_{n\in\N}$ and
$\{\xi_n\}_{n\in\N}$ decrease to zero.

We construct the sequence in order to satisfy the following list
of properties.

\begin{enumerate}[\rm (A{$_n$})]
\item $\overline{\intc P_{n-1}^{\xi_{n-1}}} \subset \intc
P_{n-1}^{\ep_{n}}\subset \overline{\intc P_{n-1}^{\ep_{n}}}\subset
\intc P_{n}^{\xi_{n}}\subset \overline{\intc P_{n}^{\xi_{n}}}
\subset \intc P_n\subset \overline{\intc P_n}\subset \intc
P_{n-1}.$

\item $1/n<\dist_{\big(\overline{\intc
P_n^{\xi_n}}\,,\,d\sigma_{\psi_n}^2\big)}(P_{n-1}^{\xi_{n-1}},P_n^{\xi_n}),$
where $d\sigma_{X_n}^2$ is the lift metric of the immersion
$\psi_n.$

\item $\psi_n(\overline{\intc P_n})\subset \B(r_n),$ where
$r_n=\sqrt{r_{n-1}^2-(2/n)^2}-\ep_n.$ Notice that \eqref{ecu: r'}
guarantees that $\{r_n\}_{n\in\N}$ decreases to a real number
$r_\infty> 1.$

\item $\|\psi_n-\psi_{n-1}\|_0<\ep_n$ in $\overline{\intc
P_{n-1}^{\ep_n}}.$

\item $\l^0_{\psi_n}\geq \a_n\cdot \l^0_{\psi_{n-1}}$ in $\overline{\intc
P_{n-1}^{\xi_{n-1}}}.$
\end{enumerate}

The sequence $\{\Upsilon_n\}_{n\in\N}$ is constructed in a
recursive way. The existence of a family $\Upsilon_1$ satisfying
assertion (C$_1$) is straightforward. The rest of the properties
have no sense for $n=1.$

Suppose that we have $\Upsilon_1,\ldots,\Upsilon_n.$ We are going
to construct $\Upsilon_{n+1}.$ We choose a decreasing sequence of
positive reals $\{\ve_m\}_{m\in\N}\searrow 0$ with
$\ve_m<\min\{1/(n+1)^2,\ep_n\}$ for all $m\in\N.$ For each $m,$ we
consider the polygon $Q_m$ and the conformal maximal immersion
$Y_m:\overline{\intc Q_m}\to\L^3$ given by Lemma \ref{lema} for
the following data:
\[
r=r_n\;,\quad P=P_n\;,\quad X=X_n\;,\quad \ep=\ve_m\;,\quad
s=\frac1{n+1}\;.
\]
For a large enough $m,$ (L.1) in Lemma \ref{lema} guarantees that
$\overline{\intc P_n^{\xi_n}}\subset \intc Q_m.$ Moreover, from
Property (L.4), we deduce that the sequence $\{Y_m\}_{m\in\N}$
uniformly converges to $\psi_n$ in $\overline{\intc
P_n^{\ve_m}}\supset \overline{\intc P_n^{\xi_n}}.$ Then, taking
into account that $Y_m$ is a harmonic map and that its Weierstrass
data are given by its derivatives, we conclude that the sequence
$\{\l^0_{Y_m}\}_{m\in\N}$ uniformly converges to $\l^0_{\psi_n}$
in $\overline{\intc P_n^{\xi_n}}.$ Hence, there exists $m_0\in\N$
satisfying
\begin{equation}\label{ecu: duke28}
\overline{\intc P_n^{\xi_n}}\subset \intc P_n^{\ve_{m_0}}\subset
\overline{\intc P_n^{\ve_{m_0}}} \subset \intc Q_{m_0}\;,
\end{equation}
\begin{equation}\label{ecu: duke29}
\l^0_{Y_{m_0}}\geq \a_{n+1}\cdot \l^0_{\psi_n}\;,\quad \text{in
}\overline{\intc P_n^{\xi_n}}\;.
\end{equation}
In order to obtain \eqref{ecu: duke29} we have taken into account
that the immersion $\psi_n$ has no branch points, it only has
singularities of lightlike type (see Remark \ref{rem: branch}).

At this point, we define $P_{n+1}=Q_{m_0},$ $\psi_{n+1}=Y_{m_0}$
and $\ep_{n+1}=\ve_{m_0}.$ From (L.2) in Lemma \ref{lema}, we
conclude that $1/(n+1)<\dist_{\big(\overline{\intc
P_{n+1}}\,,\,d\sigma_{\psi_{n+1}}^2\big)}(P_n^{\ep_{n+1}},P_{n+1}).$
Therefore, taking into account \eqref{ecu: duke28} we can take
$\xi_{n+1}$ small enough so that (A$_{n+1}$) and (B$_{n+1}$) hold.
Properties (C$_{n+1}$) and (D$_{n+1}$) are consequence of (L.3)
and (L.4), respectively, whereas \eqref{ecu: duke29} implies
(E$_{n+1}$). This concludes the construction of the sequence
$\{\Upsilon_n\}_{n\in\N}.$
\\

Now, define $\Delta:=\cup_{n\in\N}\intc P_n^{\ep_{n+1}}=
\cup_{n\in \N}\intc P_n^{\xi_n}.$ Since (A$_n$), the set $\Delta$
is an expansive union of simply connected domains resulting in
$\Delta$ being simply connected. Moreover, $\Delta$ is bounded
since Properties (A$_n$), $n\in\N$, so it is biholomorphic to a
disk. On the other hand, from (D$_n$) we obtain that
$\{\psi_n\}_{n\in\N}$ is a Cauchy sequence, uniformly on compact
sets of $\Delta.$ Then, Harnack's Theorem guarantees the existence
of a harmonic map $\psi:\Delta\to\L^3$ such that
$\{\psi_n\}_{n\in\N}\to\psi,$ uniformly on compact sets of
$\Delta.$ Then, $\psi$ has the following properties.
\\

$\bullet$ {\em $\psi$ is maximal and conformal.} This facts are
consequence of that $\psi$ is harmonic.
\\

$\bullet$ {\em $\psi$ has no branch points.} For any $z\in\Delta$
there exists $n\in\N$ so that $z\in \intc P_n^{\xi_n}.$ Given
$k>n$ and using (E$_j$), $j=n+1,\ldots,k,$ one has
$\l^0_{\psi_k}(z)\geq \a_k\cdots\a_1 \l^0_{\psi_n}(z).$ Hence,
taking the limit as $k\to\infty,$ we infer that
\[
\l^0_\psi(z)\geq \frac12 \l^0_{\psi_n}(z)>0\;,
\]
and so, $\psi$ has no branch points. Notice that the last
inequality holds because of $\psi_n$ has no branch points.

\begin{remark}
Observe that this argument does not work if we use the conformal
coefficients $\l_{\psi_k}$ instead of $\l^0_{\psi_k}.$ This fact
is implied by the possible existence of singularities of lightlike
type.
\end{remark}

$\bullet$ {\em $\psi$ is weakly complete.} This fact follows from
Properties (B$_n$), (E$_n$), $n\in\N,$ and the fact that the sum
$\sum_{n=1}^\infty 1/n$ diverges.
\\

$\bullet$ {\em $\psi(\Delta)\subset \B(1).$} Let $z\in\Delta$ and
$n\in\N$ such that $z\in \intc P_n^{\xi_n}.$ For each $k\geq n,$
Property (C$_k$) guarantees that $\psi_k(z)\in \B(r_k)\subset
\B(r_\infty).$ Taking limit as $k\to\infty,$ we obtain $\psi(z)\in
\overline{\B(r_\infty)}\subset \B(1).$
\\

This completes the proof of Theorem \ref{teorema}.


\section{Proof of Lemma \ref{lema}}\label{sec: lema}

The first step of the proof consists of the construction of a
labyrinth on $\intc P$ which depends on the polygon $P$ and a
positive integer $N.$ Let $\ell$ be the number of sides of $P.$
From now on, $N$ is a positive multiple of $\ell.$ Although $N$ is
fix, we will assume along the proof of the lemma that we have
taken it large enough so that some inequalities hold. Without loss
of generality, we assume $0\in \intc P^\ep.$

\begin{remark}
Throughout the proof of the lemma, a set of positive real
constants depending on the data of the lemma, i.e., $r,$ $P,$ $X,$
$\ep$ and $s,$ will appear. The symbol ``$\cte$'' will denote
these different constants. It is important to note that the choice
of these constants does not depend on $N.$
\end{remark}

First of all, consider $\zeta_0\in]0,\ep[.$ Therefore,
$P^{\zeta_0}$ is well defined and $\overline{\intc P^\ep}\subset
\intc P^{\zeta_0}.$ We also assume that $N$ satisfies
$2/N<\zeta_0.$

Let $v_1,\ldots,v_{2N}$ be a set of points in the polygon $P$
(containing the vertices of $P$) which divides each side of $P$
into $2N/\ell$ equal parts. Let $v'_1,\ldots, v'_{2N}$ the points
resulting from transfering the above partition to the polygon
$P^{2/N}.$ Then, we define the following sets.
\begin{itemize}
\item $L_i$ is the segment that joins $v_i$ and $v_i',$
$i=1,\ldots,2N.$

\item $\mathcal{G}_i=P^{i/N^3},$ $i=0,\ldots,2N^2.$

\item $\mathcal{A}=\bigcup_{i=0}^{N^2-1}\overline{(\intc \mathcal{G}_{2i})\setminus (\intc
\mathcal{G}_{2i+1})}$ and
$\widetilde{\mathcal{A}}=\bigcup_{i=1}^{N^2}\overline{(\intc
\mathcal{G}_{2i-1})\setminus (\intc \mathcal{G}_{2i})}.$

\item $\mathcal{R}=\bigcup_{i=0}^{2N^2}\mathcal{G}_i.$

\item $\mathcal{B}=\bigcup_{i=1}^N L_{2i}$ and $\widetilde{\mathcal{B}}=\bigcup_{i=0}^{N-1} L_{2i+1}.$

\item $\mathcal{L}=\mathcal{B}\cap \mathcal{A},$ $\widetilde{\mathcal{L}}=\widetilde{\mathcal{B}}\cap
\widetilde{\mathcal{A}}$ and
$H=\mathcal{R}\cup\mathcal{L}\cup\widetilde{\mathcal{L}}.$

\item $\Om_N=\{z\in (\intc \mathcal{G}_0)\setminus(\intc
\mathcal{G}_{2N^2})\;|\;
\dist_{(\C,\escpro{\cdot,\cdot})}(z,H)\geq 1/(4N^3)\}.$

\item $\om_i$ is the union of the segment $L_i$ and those
connected components of $\Om_N$ that have nonempty intersection
with $L_i,$ for $i=1,\ldots,2N.$

\item $\varpi_i=\{z\in\C\;|\;
\dist_{(\C,\escpro{\cdot,\cdot})}(z,\om_i)<\de(N)\},$ is chosen so
that the sets $\overline{\varpi_i},$ $i=1,\ldots,2N,$ are pairwise
disjoint.
\end{itemize}


\begin{figure}[htbp]
    \begin{center}
        \includegraphics[width=0.70\textwidth]{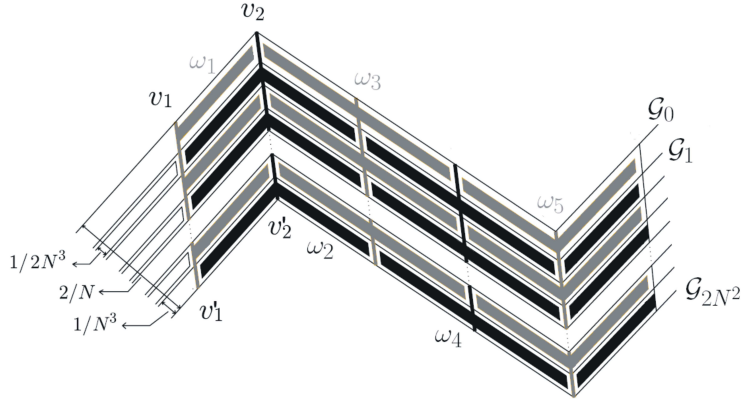}
    \end{center}
    \caption{The labyrinth.} \label{fig: 1}
\end{figure}

After constructing the labyrinth, we are going to list some of its
properties.
\begin{claim}\label{cla: abcd}
If $N$ is large enough, for any $i=1,\ldots,2N,$ one has
\begin{enumerate}[\rm A.]
\item $\diam_{(\C,\escpro{\cdot,\cdot})}(\varpi_i)<\cte/N.$

\item $\diam_{\H_+^2}
(\NN(X(\varpi_i)))<1/\sqrt{N},$ where $\diam_{\H^2_+}$ is the
intrinsic diameter in $\H^2_+.$ Here, $\NN$ is the map defined in
page \pageref{pag: normales}.

\item Denote by $(g,\Phi_3)$ the Weierstrass data of the immersion $X.$ Then,
there exists a subset $I_0\subset \{1,\ldots,2N\}$ such that
\begin{itemize}
\item $|g(z)|\neq 1$ $\forall z\in \varpi_j,$ $\forall j\in I_0.$

\item $g(z)\neq \infty$ $\forall z\in\varpi_j,$ $\forall j\in
J_0=\{1,\ldots,2N\}\setminus I_0.$
\end{itemize}

\item Let $\l^2 \escpro{\cdot,\cdot}$ be a conformal metric in
$\overline{\intc P}.$ Assume there exists $c\in\R^+$ so that
\[
\l\geq
\begin{cases}
c & \text{in }\intc P\;,
\\
c\,N^4 & \text{in }\Om_N\;.
\end{cases}
\]
Then, for any curve $\a$ in $\overline{\intc P}$ connecting
$P^{\zeta_0}$ and $P,$ one has
$\length(\a,\l\escpro{\cdot,\cdot})>\cte \,c\,N,$ where $\cte$
does not depend on $c.$
\end{enumerate}
\end{claim}
\begin{proof}
Checking Item A in the above claim is straightforward. Item B is a
consequence of Item A and the fact that $\NN$ is a differentiable
map. For a suffienctly large $N,$ Item C holds since Item A and
because of $g$ is a meromorphic function. In order to prove Item
D, we denote by $\a_j$ as the piece of $\a$ connecting $P^{j/N}$
and $P^{(j+1)/N},$ for $j=0,\ldots,N^2-1.$ Then, either the
Euclidean length of $\a_j$ is greater than $\cte/N$ or the length
of $\a_j\cap \Om_N$ is greater than $1/2N^3.$ This fact and our
assumption about $\l$ imply Item D.
\end{proof}

At this point, we construct a sequence $F_0=X,F_1,\ldots, F_{2N}$
of conformal maximal immersions (with boundary and, possibly,
lightlike singularities) defined in $\overline{\intc P}.$
\begin{claim}\label{cla: 5.2}
We will construct the sequence in order to satisfy the
following list of statements, for $i=1,\ldots,2N.$
\begin{enumerate}[\rm ({a}1$_{i}$)]
\item $F_i(z)=\Real\big( \int_0^z \phi^i(u)\,du\big)+V.$
Here, $V\in \R^3$ is a fixed vector. It does not depend on $i.$

\item $\normaR{\phi^i-\phi^{i-1}}\leq 1/N^2$ in $\overline{\intc
P}\setminus \varpi_i.$

\item $\normaR{\phi^i}\geq N^{7/2}$ in $\om_i.$

\item $\normaR{\phi^i}\geq \cte/\sqrt{N}$ in $\varpi_i.$

\item Assume $(g_i,\phi^i_3)$ are the Weierstrass data of
$F_i.$ Then, the following two assertions hold.
\begin{enumerate}[\rm ({a5.}1$_{i}$)]
\item $|g_i(z)|\neq 1,$ $\forall z\in\varpi_j,$ $\forall j\in
I_0,$ $j>i.$ Hence, the Gauss map $G_i$ of the immersion $F_i$ is
well defined in $\varpi_j$ for those $j.$ Moreover,
$\dist_{\H^2}(G_i(z),G_{i-1}(z))<1/N^2,$ for any $z\in\varpi_j$
and for any $j\in I_0,$ $j>i,$ where by $\dist_{\H^2}$ we mean the
intrinsic distance in $\H^2.$

\item $g_i(z)\neq \infty,$ $\forall z\in\varpi_j,$ $\forall j\in
J_0,$ $j>i.$ Furthermore, one has $|g_i(z)-g_{i-1}(z)|<1/N^2,$ for
any $z\in\varpi_j,$ for any of those $j.$
\end{enumerate}
\item There exists $S_i=\{e_1,e_2,e_3\}$ an orthonormal frame in
$\L^3,$ such that
\begin{enumerate}[\rm ({a6.}1$_{i}$)]
\item $\dist_{\H^2_+}(e_3,\NN(X(z))<\cte/\sqrt{N},$ for any $z\in
\overline{\varpi_i}.$

\item $(F_i(z))_{(3,S_i)}=(F_{i-1}(z))_{(3,S_i)},$ for all $z$ in
$\overline{\intc P}.$
\end{enumerate}
\item $\normaR{F_i-F_{i-1}}<\cte/N^2$ in $(\intc P)\setminus
\varpi_i.$
\end{enumerate}
\end{claim}
\begin{proof}
The sequence $F_0,F_1,\ldots,F_{2N}$ is constructed in a recursive
way. Assume that we already have $F_0,F_1,\ldots,F_{j-1}$
satisfying the assertions (a1$_i$),$\ldots,$ (a7$_i$),
$i=1,\ldots,j-1.$ Before constructing $F_j,$ we need to check the
following claim.
\begin{claim}\label{cla: 5.3}
For a large enough $N,$ the following statements hold.
\begin{enumerate}[\rm ({b}1)]
\item $\normaR{\phi^{j-1}}\leq \cte$ in $(\intc P)\setminus
(\bigcup_{k=1}^{j-1} \varpi_k).$

\item $\normaR{\phi^{j-1}}\geq \cte$ in $(\intc P)\setminus
(\bigcup_{k=1}^{j-1} \varpi_k).$

\item The diameter in $\R^3$ of $F_{j-1}(\varpi_j)$ is less than
$1/\sqrt{N}.$

\item Assume $j\in I_0.$ Then,
\begin{enumerate}
\item[\rm ({b4.}1)] The diameter in $\H^2$ of $G_{j-1}(\varpi_j)$
is less than $1/\sqrt{N}.$ In particular, there exists $p\in
G_{j-1}(\varpi_j)$ such that
$\dist_{\H^2}(p,G_{j-1}(z))<1/\sqrt{N},$ for any $z\in \varpi_j.$
\end{enumerate}
On the other hand, suppose $j\in J_0.$
\begin{enumerate}
\item[\rm ({b4.}2)] Consider the set
\[
\Gamma:=\left\{ \frac{G_{j-1}(z)}{\normaR{G_{j-1}(z)}}\;\big|\;
z\in\varpi_j\;,\;|g_{j-1}(z)|\neq 1 \right\}\;.
\]
Denote by $\Gamma^+$ (resp. $\Gamma^-$) as the part of $\Gamma$
corresponding to $\H^2_+$ (resp. $\H^2_-$). Then, there exists
$p\in\Gamma^+$ so that $\dist_{\esf^2}(\pm p,q)<1/\sqrt{N},$ for
all $q\in \Gamma^\pm.$
\end{enumerate}
\item There exists an orthonormal frame $S_j=\{e_1,e_2,e_3\}$ in
$\L^3,$ where $e_3\in \H^2_+$ and the following assertions hold.
\begin{enumerate}[\rm ({b5}.1)]
\item $\dist_{\H^2_+}(e_3,\NN(X(z)))\leq
\cte/\sqrt{N},$ for all $z\in \varpi_j.$

\item $\dist_{\H^2_+}(e_3,\pm q)\geq \cte/\sqrt{N}$ and $\dist_{\H^2_-}(-e_3,\pm q)\geq
\cte/\sqrt{N},$ for any $q$ in the set
$\{G_{j-1}(z)\;|\;z\in\varpi_j,\;|g_{j-1}(z)|\neq 1\}.$ We mean
that we only have to compute the distance if both points are in
the same connected component of $\H^2.$
\end{enumerate}
\end{enumerate}
\end{claim}
\begin{proof}
To deduce (b1) and (b2) we have to use just (a2$_k$),
$k=1,\ldots,j-1.$ Item (b3) is a consequence of (b1) and Claim
\ref{cla: abcd}.A. In order to prove (b4) we distinguish cases. If
$j\in I_0,$ taking into account Claim \ref{cla: abcd}.A and Claim
\ref{cla: abcd}.C we obtain that the diameter of $G_0(\varpi_j)$
is bounded by $\cte/N.$ Then, we can apply (a5.1$_k$),
$k=1,\ldots,j-1,$ to conclude (b4.1). On the other hand, if $j\in
J_0,$ we use again Claim \ref{cla: abcd}.A and Claim \ref{cla:
abcd}.C to deduce that $\diam_\C(g_0(\varpi_j))<\cte/N.$
Therefore, (a5.2$_k$), $k=1,\ldots,j-1,$ imply that
$\diam_\C(g_{j-1}(\varpi_j))<\cte/N.$ This fact guarantees (b4.2)
for a large enough $N.$ We also have taken into account that if
$|g_{i-1}(z)|<1<|g_{i-1}(z')|$ and $g_{i-1}(z)\approx
g_{i-1}(z'),$ then $G_{j-1}(z) \approx -G_{j-1}(z').$

The proof of (b5) is slightly more complicated. First, assume that
$j\in I_0.$ Without loss of generality we can assume that
$G_{j-1}(\varpi_j)\subset \H_+^2,$ otherwise we would work with
$-G_{j-1}(\varpi_j).$ Consider $p$ given by Property (b4.1), then
to obtain (b5.2), it suffices to take $e_3$ in
$C=\{q\in\H^2_+\;|\;\dist_{\H^2_+}(p,q)>2/\sqrt{N}\}.$ Moreover,
in order to satisfy (b5.1), the vector $e_3$ must be chosen as
follows.
\begin{itemize}
\item If $C\cap \NN(X(\varpi_j))\neq \emptyset,$ then we take
$e_3$ in that set. Therefore (b5.1) holds because of Claim
\ref{cla: abcd}.B.

\item If $C\cap \NN(X(\varpi_j))= \emptyset,$ then we take
$e_3\in C$ such that $\dist_{\H^2_+}(e_3,q')<2/\sqrt{N}$ for some
$q'\in \NN(X(\varpi_j)).$ This choice is possible since (b4.1).
Again Claim \ref{cla: abcd}.B. guarantees (b5.1).
\end{itemize}

Assume now that $j\in J_0.$ We define the sets
\[
\Lambda_\pm:=\left\{\frac{q}{\normaR{q}}\;\big|\;q\in\H^2_\pm\right\}\subset\esf^2\;,\quad
\Xi:=\left\{ \frac{\NN(X(z))}{\normaR{\NN(X(z))}}
\;\big|\;z\in\varpi_j\right\}\subset\Lambda_+\;.
\]
In order to prove assertion (b5) in this case, we are going to use
the following statement. There exists $e_3\in \H^2_+$ so that the
vector $\widehat{e}_3=e_3/\normaR{e_3}$ satisfies
\begin{enumerate}[ i)]
\item $\dist_{\Lambda_+}(\widehat{e}_3,q)\leq \cte/\sqrt{N},$ for all
$q\in \Xi.$

\item $\dist_{\Lambda_+}(\widehat{e}_3,\pm q)\geq \cte/\sqrt{N}$ and
$\dist_{\Lambda_-}(-\widehat{e}_3,\pm q)\geq \cte/\sqrt{N}$ for
any $q\in\Gamma.$ Again, we mean that we only have to compute the
distance if both points are in $\Lambda_+$ or both in $\Lambda_-.$
\end{enumerate}
Indeed, the proof consists of the same arguments as above but
using (b4.2) instead of (b4.1). Then, (b5.1) is a consequence of
i) and the fact that $\normaR{\NN(X(\varpi_j))}$ is bounded (not
depending on $N$). Moreover, ii) implies (b5.2). Hence, $e_3$
proves Property (b5) in this case.
\end{proof}

Now, we can continue with the proof of Claim \ref{cla: 5.2}. Let
$(g^{j-1},\phi_3^{j-1})$ be the Weierstrass data of the immersion
$F_{j-1}$ in the basis $S_j$ given by (b5). For any $\a>0,$
consider $h_\a:\overline{\intc P}\to\C$ a holomorphic function
without zeros and satisfying
\begin{itemize}
\item $|h_\a-1|<1/\a$ in $\overline{\intc P}\setminus \varpi_j.$

\item $|h_\a-\a|<1/\a$ in $\om_j.$
\end{itemize}
This family of functions is given by Runge's Theorem. Using $h_\a$
as a López-Ros parameter, we define $F_j$ in the coordinate system
$S_j$ as $g^j=g^{j-1}/h_\a$ and $\phi_3^j=\phi_3^{j-1}.$ Taking
into account that $h_\a\to 1$ (resp. $h_\a\to\infty$) uniformly in
$\overline{\intc P}\setminus \varpi_j$ (resp. in $\om_j$), as
$\a\to\infty,$ it is clear that Properties (a1$_j$), (a2$_j$),
(a3$_j$), (a5$_j$) and (a7$_j$) hold for a large enough (in terms
of $N$) value of the parameter $\a.$ Moreover, using (b5.1) we
obtain (a6.1$_j$) and to get (a6.2$_j$) we use that
$\phi_3^{j-1}=\phi_3^j$ in the frame $S_j.$ Finally, we are going
to prove (a4$_j$). Consider $z\in\varpi_j$ with $|g^{j-1}(z)|\neq
1.$ Using the stereographic projection for $\H^2$ from the point
$e_3\in\H^2_+,$ from Property (b5.2) one has
\[
\frac{\sinh \big(\frac{\cte}{\sqrt{N}}\big)}{\cosh
\big(\frac{\cte}{\sqrt{N}}\big)+1}\leq |g^{j-1}(z)|\leq
\frac{\sinh \big(\frac{\cte}{\sqrt{N}}\big)}{\cosh
\big(\frac{\cte}{\sqrt{N}}\big)-1}\;.
\]
On the other hand, if $|g^{j-1}(z)|=1,$ then the above
inequalities trivially hold, so they occur for any $z\in\varpi_j.$
Therefore,
\[
\normaR{\phi^j} \geq |\phi_3^j|=|\phi_3^{j-1}|\geq \sqrt{2}
\normaR{\phi^{j-1}}\frac{|g^{j-1}|}{1+|g^{j-1}|^2}
\]
\[
\geq \cte \cdot \tanh \left(\frac{\cte}{\sqrt{N}}\right)\geq
\frac{\cte}{\sqrt{N}}\quad \text{in }\varpi_j\;,
\]
where we have used (a6.2$_j$) and (b2). This fact proves (a4$_j$)
and concludes the proof of Claim \ref{cla: 5.2}.
\end{proof}

\begin{remark}
Notice that in the definition of $F_i$ in Property (a1$_i$), we
need the addition of the fixed vector $V.$ Otherwise, it would be
$F_i(0)=(0,0,0).$ In particular, $X(0)=(0,0,0)\notin \B(r),$ which
is absurd.
\end{remark}

\begin{remark}\label{rem: ortonormal}
Let $S_i=\{e_1,e_2,e_3\}$ be the $\L^3$-orthonormal basis given by
Property (a6$_i$). Consider
$\widetilde{S}_i=\{\widetilde{e}_1,\widetilde{e}_2,\widetilde{e}_3\}$
an $\R^3$-orthonormal basis such that $\{e_1,e_2\}$ and
$\{\widetilde{e}_1,\widetilde{e}_2\}$ define the same plane, and
$e_3$ and $\widetilde{e}_3$ lie in the same halfspace determined
by that plane, i.e.,
$\widetilde{e}_3=-\mathcal{J}(e_3)/\normaR{e_3},$ where
$\mathcal{J}(e_3^1,e_3^2,e_3^3)=(e_3^1,e_3^2,-e_3^3).$ Then, one
has
\begin{itemize}
\item
$\dist_{\esf^2}(\widetilde{e}_3,\NN_0(X(z)))<\cte/\sqrt{N},$ for
any $z\in\varpi_i,$ where $\NN_0$ is the map that was defined in
page \pageref{pag: normales}.

\item $(F_i(z))_{(3,\widetilde{S}_i)}=
(F_{i-1}(z))_{(3,\widetilde{S}_i)}.$
\end{itemize}
\end{remark}

Now, we establish some properties of the final immersion $F_{2N}.$

\begin{claim}\label{cla: final}
If $N$ is large enough, then $F_{2N}$ satisfies
\begin{enumerate}[\rm ({c}1)]
\item $2s < \dist_{(\overline{\intc
P}\,,\,d\sigma_{F_{2N}})}(P,P^\ep),$ where by $d\sigma_{F_{2N}}$
we represent the lift metric of the immersion $F_{2N}.$

\item $\normaR{F_{2N}-X}<\cte/N,$ in $\overline{\intc P}\setminus
(\cup_{i=1}^{2N} \varpi_i).$

\item There exists a polygon $Q$ such that
\begin{enumerate}[\rm ({c3}.1)]
\item $\overline{\intc P^\ep} \subset \intc Q \subset
\overline{\intc Q}\subset \intc P.$

\item $s<\dist_{(\overline{\intc P}\,,\,
d\sigma_{F_{2N}})}(z,P^\ep)<2s,$ for any $z\in Q.$

\item $F_{2N}(\overline{\intc Q})\subset \B(R),$
where $R=\sqrt{r^2-4s^2}-\ep.$
\end{enumerate}
\end{enumerate}
\end{claim}
\begin{proof}
Properties (b2), (a2$_i$), (a3$_i$) and (a4$_i$), $i=1,\ldots,2N,$
guarantee that the conformal coefficient $\l_{F_{2N}}^0$ of the
lift metric of $F_{2N}$ satisfies
\[
\l_{F_{2N}}^0 =\frac{\normaR{\phi^{2N}}}{\sqrt{2}}\geq
\begin{cases}
\frac{\cte}{\sqrt{N}} & \text{in }\intc P\\
\frac{\cte}{\sqrt{N}}\, N^4 & \text{in }\Omega_N\;.
\end{cases}
\]
Therefore, Claim \ref{cla: abcd}.D imply that
\[
\dist_{(\overline{\intc P}\,,\,d\sigma_{F_{2N}})}(P,P^\ep)\geq
\dist_{(\overline{\intc P}\,,\,d\sigma_{F_{2N}})}(P,P^{\zeta_0})>
\frac{\cte}{\sqrt{N}}\,N=\cte\sqrt{N}>2s\;,
\]
for a large enough $N.$ We have proved (c1). Property (c2)
trivially holds from (a2$_i$), $i=1,\ldots,2N.$

In order to construct the polygon $Q$ of the assertion (c3), we
consider the set
\[
\mathcal{K}=\big\{ z\in (\intc P)\setminus (\intc P^\ep)\;\big|\;
s< \dist_{(\overline{\intc P}\,,\,d\sigma_{F_{2N}})}(z,P^\ep)<2s
\big\}\;.
\]
From (c1), $\mathcal{K}$ is a nonempty open subset of $(\intc
P)\setminus (\intc P^\ep),$ and $P$ and $P^\ep$ are contained in
different connected components of $\C\setminus\mathcal{K}.$
Therefore, we can choose a polygon $Q$ on $\mathcal{K}$ satisfying
(c3.1) and (c3.2).

The proof of (c3.3) is more complicated. Consider
$z\in\overline{\intc Q}.$ First, we assume that $z\in (\intc
P)\setminus (\cup_{i=1}^{2N}\varpi_i).$ Then, we can use
Properties (a2$_i$), $i=1,\ldots,2N,$ to conclude that
$\normaR{F_{2N}(z)-X(z)}<\cte/N.$ Moreover, from the hypotheses of
Lemma \ref{lema}, we have $X(z)\in \B(r).$ Hence, $F_{2N}(z)\in
\B(R),$ if $N$ is large enough.

On the other hand, suppose that there exists $i\in
\{1,\ldots,2N\}$ with $z\in\varpi_i.$ Choose a curve $\g:[0,1]\to
\intc P$ satisfying $\g(0)\in P^\ep,$ $\g(1)=z$ and
$\length(\g,d\sigma_{F_{2N}})< 2s.$ This election is possible
since (c3.2). Label
\[
t_0=\sup\big\{ t\in[0,1]\;\big|\; \g(t)\in\partial
\varpi_i\big\}\;,\quad z_0=\g(t_0)\;.
\]
Notice that this supremum exists because $\varpi_i\subset (\intc
P)\setminus \overline{\intc P^\ep}$ (for a large enough $N$). Now,
consider the basis $\widetilde{S}_i$ explained in Remark \ref{rem:
ortonormal}, then we have
\begin{equation}\label{ecu: pita12}
\norma{(F_{2N}(z)-X(z))_{(*,\widetilde{S}_i)}}\leq
2s+\frac{\cte}{\sqrt{N}}\;,
\end{equation}
\begin{equation}\label{ecu: pita3}
|(F_{2N}(z)-X(z))_{(3,\widetilde{S}_i)}|<\frac{\cte}N\;.
\end{equation}
\begin{figure}[htbp]
    \begin{center}
        \includegraphics[width=0.50\textwidth]{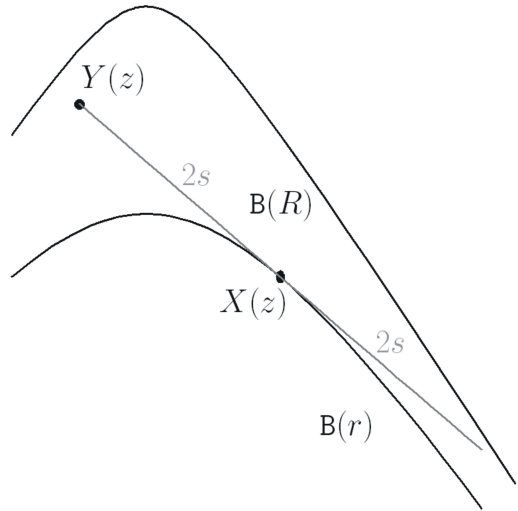}
    \end{center}
    \caption{The effect of the deformation.}\label{fig: 2}
\end{figure}
Indeed,
\[
\norma{(F_{2N}(z)-X(z))_{(*,\widetilde{S}_i)}}\leq\normaR{F_{2N}(z)-F_{2N}(z_0)}+
\normaR{F_{2N}(z_0)-F_{i-1}(z_0)}+
\]
\[
\normaR{F_{i-1}(z_0)-F_{i-1}(z)}+\normaR{F_{i-1}(z)-X(z)}\leq
\]
\[
\length(\g,d\sigma_{F_{2N}})+\frac{\cte}N+\frac1{\sqrt{N}}+\frac{\cte}N<2s+\frac{\cte}{\sqrt{N}}\;,
\]
where we have used (a7$_j$), $j=1,\ldots,2N,$ and (b3). On the
other hand, taking Remark \ref{rem: ortonormal} and (a7$_j$),
$j=1,\ldots,2N,$ into account, we conclude
\[
|(F_{2N}(z)-X(z))_{(3,\widetilde{S}_i)}|\leq
\normaR{F_{2N}(z)-F_i(z)}+
|(F_{i}(z)-F_{i-1}(z))_{(3,\widetilde{S}_i)}|+
\]
\[
\normaR{F_{i-1}(z)-X(z)}<\frac{\cte}N+\frac{\cte}N=\frac{\cte}N\;.
\]

At this point, consider the following statement. Its proof is
elemental, we leave the details to the reader.
\begin{claim}\label{cla: bolas}
Let $0<x<t.$ Consider $p\in \B(t)$ and $v\in\R^3$ with
$\escproR{\NN_0(p),v}=0$ and $\normaR{v}=x.$ Then, $p+v\in
\B(\sqrt{t^2-x^2}).$
\end{claim}

Now, Remark \ref{rem: ortonormal}, equations \eqref{ecu: X-cota},
\eqref{ecu: pita12} and \eqref{ecu: pita3}, and the above claim
guarantee that $F_{2N}(z)\in \B(R),$ if $N$ was chosen large
enough. This proves (c3.3) and finishes the proof of Claim
\ref{cla: final}.
\end{proof}

From Claim \ref{cla: final} it is straightforward to check that
(for $N$ large enough) $Y=F_{2N}:\overline{\intc Q}\to\L^3$ proves
Lemma \ref{lema}.


\end{document}